\documentclass[11pt,reqno]{article}
\usepackage{amsmath,mathrsfs,graphicx,amssymb,amsfonts,color}

\font\bg=cmbx10 scaled\magstep1

\font\small=cmr8

\newtheorem{newlemma}{{\bf Lemma}}

\newenvironment{lema}{\begin{newlemma}{\hspace{-0.5
em}{\bf.}}}{\end{newlemma}}

\newtheorem{newteorem}{{\bf Theorem}}

\newenvironment{teorem}{\begin{newteorem}{\hspace{-0.5
em}{\bf.}}}{\end{newteorem}}

\newtheorem{newkorolari}{{\bf Corollary}}

\newenvironment{korolari}{\begin{newkorolari}{\hspace{-0.5
em}{\bf.}}}{\end{newkorolari}}

\newtheorem{newdefine}{{\bf Definition}}

\newtheorem{newquestion}{{\bf Question}}

\newtheorem{newkonjek}{{\bf Conjecture}}

\newtheorem{newexample}{{\bf Example}}

\begin{document}
\tolerance=10000
\baselineskip18truept
\newbox\thebox
\global\setbox\thebox=\vbox to 0.2truecm{\hsize 0.15truecm\noindent\hfill}
\def\boxit#1{\vbox{\hrule\hbox{\vrule\kern0pt
\vbox{\kern0pt#1\kern0pt}\kern0pt\vrule}\hrule}}
\def\qed{\lower0.1cm\hbox{\noindent \boxit{\copy\thebox}}\bigskip}
\def\ss{\smallskip}
\def\ms{\medskip}
\def\bs{\bigskip}
\def\c{\centerline}
\def\nt{\noindent}
\def\ul{\underline}
\def\ol{\overline}
\def\lc{\lceil}
\def\rc{\rceil}
\def\lf{\lfloor}
\def\rf{\rfloor}
\def\ov{\over}
\def\t{\tau}
\def\th{\theta}
\def\k{\kappa}
\def\l{\lambda}
\def\L{\Lambda}
\def\g{\gamma}
\def\d{\delta}
\def\D{\Delta}
\def\e{\epsilon}
\def\lg{\langle}
\def\rg{\rangle}
\def\p{\prime}
\def\sg{\sigma}
\def\ch{\choose}

\newcommand{\ben}{\begin{enumerate}}
\newcommand{\een}{\end{enumerate}}
\newcommand{\bit}{\begin{itemize}}
\newcommand{\eit}{\end{itemize}}
\newcommand{\bea}{\begin{eqnarray*}}
\newcommand{\eea}{\end{eqnarray*}}
\newcommand{\bear}{\begin{eqnarray}}
\newcommand{\eear}{\end{eqnarray}}

\centerline{\Large \bf  On the domination polynomials of  cactus chains }
\bigskip

\bs

\baselineskip12truept
\centerline{S. Alikhani$^{}${}\footnote{\baselineskip12truept\it\small
Corresponding author. E-mail: alikhani@yazd.ac.ir},  S. Jahari and M. Mehryar }
\baselineskip20truept
\centerline{\it Department of Mathematics, Yazd University}
\vskip-8truept
\centerline{\it  89195-741, Yazd, Iran}

\vskip-0.2truecm
\nt\rule{16cm}{0.1mm}

\nt{\bg ABSTRACT}
\medskip

\baselineskip14truept

\nt{ Let $G$ be a simple graph of order $n$.
The domination polynomial of $G$ is the polynomial
$D(G, x)=\sum_{i=\gamma(G)}^{n} d(G,i) x^{i}$,
where $d(G,i)$ is the number of dominating sets of $G$ of size $i$ and
$\gamma(G)$ is the domination number of $G$. In this paper we consider cactus chains with triangular and
square blocks and study their domination polynomials.}

\ms

\nt{\bf Mathematics Subject Classification:} {\small 05C60, 05C69.}

\nt\rule{16cm}{0.1mm}

\baselineskip20truept

\section{Introduction}

\nt Let $G=(V,E)$ be a simple graph.
For any vertex $v\in V(G)$, the {\it open neighborhood} of $v$ is the
set $N(v)=\{u \in V (G) | \{u, v\}\in E(G)\}$ and the {\it closed neighborhood} of $v$
is the set $N[v]=N(v)\cup \{v\}$. For a set $S\subseteq V(G)$, the open
neighborhood of $S$ is $N(S)=\bigcup_{v\in S} N(v)$ and the closed neighborhood of $S$
is $N[S]=N(S)\cup S$.
A set $S\subseteq V(G)$ is a {\it dominating set} if $N[S]=V$ or equivalently,
every vertex in $V(G)\backslash S$ is adjacent to at least one vertex in $S$.
The {\it domination number} $\gamma(G)$ is the minimum cardinality of a dominating set in $G$.
For a detailed treatment of these parameters, the reader is referred to~\cite{domination}.
Let ${\cal D}(G,i)$ be the family of dominating sets of a graph $G$ with cardinality $i$ and
let $d(G,i)=|{\cal D}(G,i)|$.
The {\it domination polynomial} $D(G,x)$ of $G$ is defined as
$D(G,x)=\sum_{ i=\gamma(G)}^{|V(G)|} d(G,i) x^{i}$,
where $\gamma(G)$ is the domination number of $G$ (see \cite{euro,saeid1}). Obviously, the number of dominating sets of a graph $G$ is
$D(G,1)$ (see \cite{gcom,gcomkot}). Recently the number of the dominating sets of graph $G$, i.e., $D(G,1)$ has been considered and studied in \cite{dam} with a different
approach.

\nt Domination theory  have many applications in sciences and technology (see \cite{domination}). Recently the dominating set  has found
application in  the assignment of structural domains in complex protein structures, which is an important task in bio-informatics (\cite{eslahchi}).


\ms

\nt We recall that the Hosoya
index $Z(G)$ of a molecule graph $G$, is the  number of matching sets,  and the Merrifield-Simmons index $i(G)$ of graph $G$, is the  number
of independent sets. The Hosoya index of a graph has application to correlations with boiling points, entropies, calculated
bond orders, as well as for coding of chemical structures. The Merrifield-Simmons index is one of the most popular topological indices in chemistry.
For more information of these two indices see \cite{mesi,prod,wag}. Note that $Z(G)$ and $i(G)$ can be study by  the value of matching polynomial and independence polynomial at $1$.

\ms

\nt In this paper we consider a  class of simple linear polymers called cactus chains. Cactus graphs were first known as Husimi
trees; they appeared in the scientific literature some sixty years ago in papers by Husimi and
Riddell concerned with cluster integrals in the theory of condensation in statistical mechanics \cite{9,12,14}.
 We refer the reader to papers \cite{chellali,13} for some aspects of domination in cactus
graphs.

\nt A cactus graph is a connected graph in which no edge lies in more than one cycle. Consequently,
each block of a cactus graph is either an edge or a cycle. If all blocks of a cactus $G$
are cycles of the same size $i$, the cactus is $i$-uniform.
A triangular cactus is a graph whose blocks are triangles, i.e., a $3$-uniform cactus.
 A vertex shared by two or more triangles is called a cut-vertex. If each triangle of a triangular
cactus $G$ has at most two cut-vertices, and each cut-vertex is shared by exactly two triangles,
we say that $G$ is a chain triangular cactus. By replacing triangles in this  definitions by cycles of length $4$ we obtain cacti whose
every block is $C_4$.
 We call such cacti square cacti. Note that the internal squares may differ in the way they connect to their
 neighbors. If their cut-vertices are adjacent, we say that such a square is an ortho-square;
if the cut-vertices are not adjacent, we call the square a para-square.

\nt In Section 2 we study the domination polynomial of the chain triangular cactus with two approach.   In Section 3 we study
the domination polynomials of chains of squares.

\section{Domination polynomials of the chain triangular cactus }

\begin{figure}[!h]
\hspace{4cm}
\includegraphics[width=7.3cm,height=1.5cm]{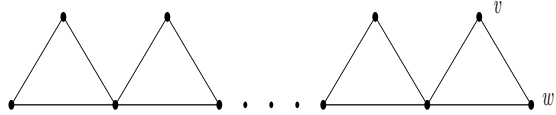}
\caption{ \label{cactus} The chain triangular cactus.}
\end{figure}

 We call the number of triangles in $G$,  the
length of the chain. An example of a chain triangular cactus is shown in Figure \ref{cactus}.
Obviously, all chain triangular cacti of the same length are isomorphic.
Hence, we denote the chain triangular cactus of length $n$ by $T_n$.
 In  this paper  we investigate the domination polynomial of $T_n$ by two different approach.

\subsection{Computation of $D(T_n,x)$ using recurrence relation}

\nt In the first subsection, we use
results and recurrence relations of the domination polynomial of a graph to find a recurrence relation for $D(T_n,x)$.

\nt We need the following theorem:
\begin{teorem}\label{theorem4}{\rm \cite{ saeid1}}
If a graph $G$ consists of $k$ components $G_1, \dots ,G_k$,  then
$D(G, x) =\prod_{i=1}^k D(G_i, x).$
\end{teorem}

\nt The vertex contraction $G/u$ of a graph $G$ by a vertex $u$ is the operation under
which all vertices in $N(u)$ are joined to each other and then $u$ is deleted (see\cite{Wal}).

\nt The following theorem is  useful for finding the recurrence relations for the  domination polynomials  of  arbitrary graphs.

\begin{teorem}\label{theorem1}{\rm \cite{saeid2,Kot}}
Let $G$ be a graph. For any vertex $u$ in $G$ we have
\begin{eqnarray*}
D(G, x) = xD(G/u, x) + D(G - u, x) + xD(G - N[u], x) - (1 + x)p_u(G, x),
\end{eqnarray*}
where $p_u(G, x)$ is the polynomial counting the dominating sets of $G - u$ which do not contain any
vertex of $N(u)$ in $G$.
\end{teorem}

\nt Domination polynomial satisfies a recurrence relation for arbitrary graphs which is based on the edge and vertex
elimination operations. The recurrence uses composite operations, e.g. $G - e /u$, which
stands for $(G - e) /u$.

\begin{teorem}\label{theorem3}{\rm \cite{ Kot}}
Let $G$ be a graph. For every edge $e =\{u, v\}\in E$,
\begin{eqnarray*}
D(G, x)&=&D(G - e, x) +\frac{x}{x-1} \Big[ D(G - e/u, x) + D(G - e/v, x)\\
&-&D(G/u, x) - D(G/v, x) - D(G - N[u], x) - D(G - N[v], x)\\
&+&D(G - e - N[u], x) + D(G - e - N[v], x)\Big].
\end{eqnarray*}
\end{teorem}


\nt  We use for graphs $G = (V, E)$ the following vertex operation,
which is commonly found in the literature. Let $v \in V$ be a vertex of $G$.
A vertex appending $G+e$ (or $G+ \{v, \cdot\}$) denotes the graph ($V \cup \{v'\}, E\cup \{v, v'\}$) obtained from
$G$ by adding a new vertex $v'$ and an edge $\{v, v'\}$ to $G$.

\nt The following theorem gives recurrence relation for the domination polynomial of $T_n$.

\begin{teorem}\label{theorem5}
For every $n\geq 3$,
$$D(T_n,x)=(x^2+2x)D(T_{n-1},x)+(x^2+x)D(T_{n-2},x),$$ with initial condition $D(T_1,x)=x^3+3x^2+3x$
and $D(T_2,x)=x^5+5x^4+10x^3+8x^2+x$.
\end{teorem}

\nt{\bf Proof.}
\nt Consider the graph $T_n$  as shown in  the following Figure ~\ref{cactus}. Since $T_n/u$ is isomorphic to $T_n - u$ and $p_u(T_n, x)=0$, by Theorem~\ref{theorem1} we have:
\begin{eqnarray}\label{eq1}
D(T_n, x) &=& x D(T_n/u, x) + D(T_n - u, x) + xD(T_n - N[u], x) - (1 + x)p_u(T_n, x)\nonumber\\
&=& (x+1) D(T_n/u, x)  + xD(T_n - N[u], x)\nonumber\\
&=& (x+1) D(T_{n-1}+ e, x)  + xD(T_{n-2}+ e, x).
\end{eqnarray}
\nt

\begin{figure}[!h]
\hspace{3.7cm}
\includegraphics[width=7cm,height=1.8cm]{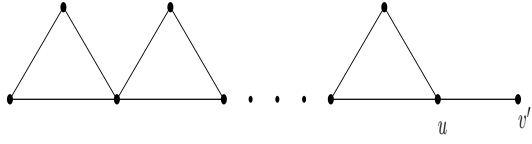}
\caption{ \label{cactus.v}      The Graph $T_{n-1}+e.$ }
\end{figure}

\nt Note  we use Theorems ~\ref{theorem4} and~\ref{theorem1} to obtain the domination polynomial of the graph
$T_{n-1}+e$ (see Figure \ref{cactus.v}). Suppose that $v'$ be a vertex of degree $1$ in  graph $T_{n-1}+e$ and let $u$ be its neighbor.  Note that in this case $p_u(T_{n-1}+e,x)=0$. We deduce that for each $n \in \mathbb{N}$,
$D(T_{n-1}+ e, x) =x [ D(T_{n-1}, x) + D(T_{n-2}+e, x)+D(T_{n-3}+e,x)]$. Therefore by equation (\ref{eq1}) and this equality we have
$$D(T_n,x)=(x^2+x)(D(T_{n-1},x)+D(T_{n-3}+e,x))+ (x^2+2x) D(T_{n-2}+e,x).$$
\nt Now it's suffices to prove the following equality:
 $$(x^2+x)D(T_{n-3}+e,x)+(x^2+2x)D(T_{n-2}+e,x)=xD(T_{n-1},x)+(x^2+x)D(T_{n-2},x).$$
 For this purpose we use  Theorem \ref{theorem1} for
$D(T_{n-1},x)$. We have
$$xD(T_{n-1},x)=(x^2+x)D(T_{n-2}+e,x)+
x^2D(T_{n-3}+e,x).$$
Now we use Theorem~\ref{theorem1} for $v'$ to obtain domination
polynomial of $T_{n-2}+e$, then we have
$D(T_{n-2}+e,x)=
(1+x)D(T_{n-2},x)+xD(T_{n-3}+e,x)-(1+x)D(T_{n-3}+e,x).$ Therefore the result follows.\quad\qed

\subsection{Computation of $D(T_n,x)$ by counting the number of dominating sets}
\nt In this section we shall obtain a recurrence relation for the domination polynomial of $T_n$. For this purpose we
count the number of dominating sets of $T_n$ with cardinality $k$.
In other words, we first find a two variables recursive formula for $d(T_n,k)$.

\nt Recently by private communication, we found that the following result also appear in \cite{doslicnew} but were proved independently.  

\begin{teorem}\label{recurrence}
The number of dominating sets of  $T_n$ with cardinality $k$ is given by
$$d(T_n,k)=2d(T_{n-1},k-1)+d(T_{n-1},k-2)+d(T_{n-2},k-1)+d(T_{n-2},k-2).$$
\end{teorem}
\nt{\bf Proof.}
We shall make a dominating set of $T_n$ with cardinality $k$ which we denote it by $\mathcal{T}_{n}^k$. We consider all cases:

\nt {\bf Case 1.} If $\mathcal{T}_{n}^k$ contains both of $v$ and $w$, then we have $\mathcal{T}_{n}^k=\mathcal{T}_{n-1}^{k-2} \cup \{v,w\}$.
In this case we have $d(T_n,k)=d(T_{n-1},k-2)$.

\nt {\bf Case 2.} If  $\mathcal{T}_{n}^k$ contains only $v$ or $w$ (say $v$), then we have $\mathcal{T}_{n}^k=\mathcal{T}_{n-1}^{k-1} \cup \{v\}$.
In this case we have $d(T_n,k)=2d(T_{n-1},k-1)$.

\begin{figure}[!h]
\hspace{2cm}
\includegraphics[width=10cm,height=8cm]{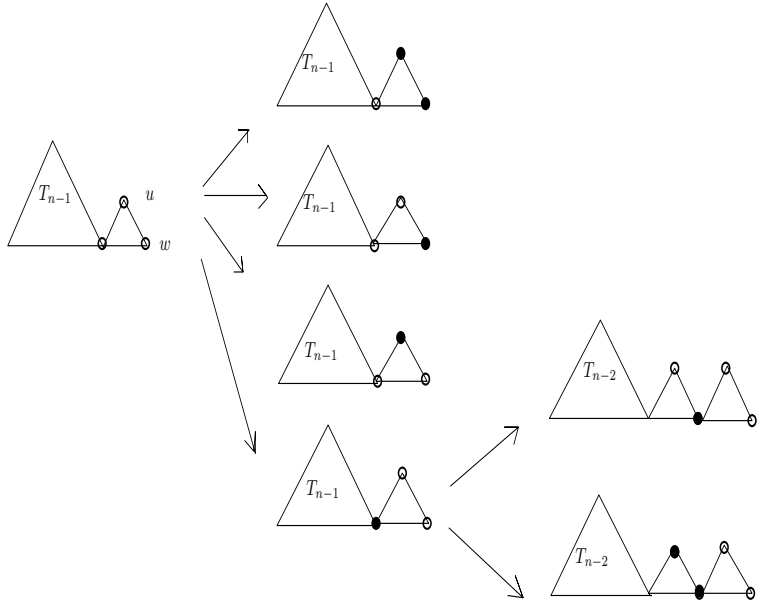}
\caption{ \label{tn} Recurrence relation for $d(T_n,k)$. }
\end{figure}

\nt {\bf Case 3.} If $\mathcal{T}_{n}^k$  contains none of $v$ and $w$, then we can construct $\mathcal{T}_{n}^k$ by
$\mathcal{T}_{n-2}^{k-1}$ or $\mathcal{T}_{n-2}^{k-2}$ as shown   in Figure \ref{tn}. In this case we have
$d(T_n,k)=d(T_{n-2},k-1)+d(T_{n-2},k-2)$.
 By adding all contributions we obtain the recurrence for $d(T_n,k)$.\quad\qed

\begin{korolari}\label{main}
For every $n\geq 3$,
$$D(T_n,x)=(x^2+2x)D(T_{n-1},x)+(x^2+x)D(T_{n-2},x).$$
\end{korolari}
\nt{\bf Proof.} It follows from Theorem \ref{recurrence} and the definition of the domination polynomial.\quad\qed

\ms

\nt We mention here the Hosoya index of a graph $G$ is the total number of matchings of $G$  and the Merrifield-Simmons index is
the total number of its independent sets. Motivation by these indices, we are interested to count the total
number of dominating set of a graph which is equal to   $D(G,1)$. Here we present a recurrence relation to the total number of the chain
triangular cactus.

\begin{teorem}
The enumerating sequence $\{t_n\}$ for the number of dominating sets in $T_n$ $(n\geq 2)$ is
$$t_n= 3t_{n-1} + 2t_{n-2}$$
with  initial values  $t_0 = 2$, $t_1 = 7$.
\end{teorem}
\nt{\bf Proof.} Since $t_n=D(T_n,1)$, it follows from Corollary \ref{main}.\quad\qed

\section{Domination polynomials of chains of squares}

\nt By replacing triangles in the definitions of triangular cactus,  by cycles of length $4$ we obtain cacti whose
every block is $C_4$. We call such cacti, square cacti. An example of a square cactus chain is
shown in Figure \ref{pcactus}. We see that the internal squares may differ in the way they connect to their neighbors. If their cut-vertices are adjacent, we say that such a square is an ortho-square;
if the cut-vertices are not adjacent, we call the square a para-square.

\subsection{Domination polynomial of para-chain square cactus graphs}

\begin{figure}[!h]
\hspace{3.9cm}
\includegraphics[width=7.3cm,height=2cm]{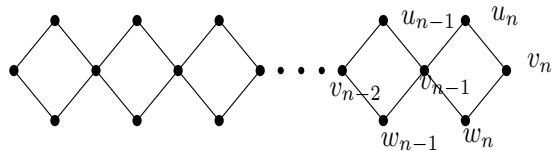}
\caption{ \label{pcactus} Para-chain square cactus graphs. }
\end{figure}

\nt In this subsection we consider a para-chain of length $n$, $Q_n$,  as shown in Figure \ref{pcactus}. We shall
obtain a recurrence relation for the domination polynomial of $Q_n$. As usual we denote the number of dominating
sets of $Q_n$ by $d(Q_n,k)$. The following theorem gives a  recurrence relation for $D(Q_n,x)$.

\nt We need the following Lemma for finding domination polynomial of the $Q_n$.

\begin{figure}[!h]
\hspace{1cm}
\includegraphics[width=14.6cm,height=1.9cm]{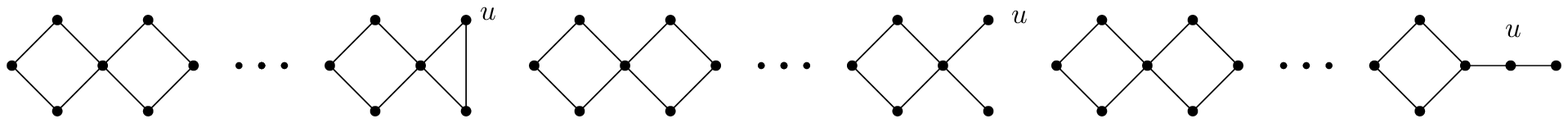}
\caption{\label{figure1} Graphs $Q_n^\bigtriangleup$, ~$Q_n'$~and ~$Q_n(2)$, respectively.}
\end{figure}

\begin{figure}[!h]
\hspace{1.3cm}
\includegraphics[width=14cm,height=1.9cm]{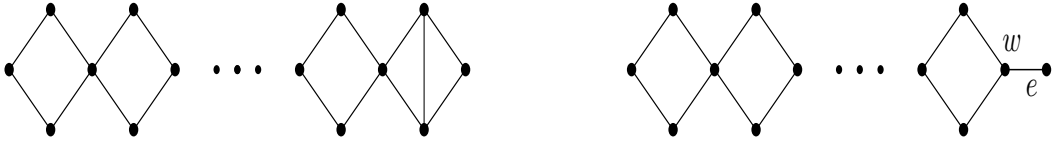}
\caption{\label{figure2} Graphs ~$(Q_n+e)/w$~and~ $Q_n+e$, respectively.}
\end{figure}

\begin{lema}\label{lemma1}
\nt For  graphs in figures~\ref{figure1} and \ref{figure2} have:

\nt $(i)~ D(Q_n^\bigtriangleup, x)=(1+x)D(Q_n+e,x) + xD(Q_{n-1}',x)$, where $D(Q_0^\bigtriangleup, x)=x^3+3x^2+3x$.

\nt $(ii)~ D(Q_n(2), x)=x(D(Q_n+e,x) + D(Q_n,x) + D(Q_{n-1}',x))$, where $D(Q_0(2), x)=x^3+3x^2+x$.

\nt $(iii)~ D(Q_n', x)=(1+x)D(Q_n+e,x) - xD(Q_{n-1}',x)$, where $D(Q_0', x)=x^3+3x^2+x$.

\nt $(iv)~ D(Q_n+e, x)=x(D(Q_n,x) + D(Q_{n-1},x)) + xD(Q_{n-1}',x) + 2x^2D(Q_{n-2}',x)$, where $D(Q_1+e, x)=x^5+5x^4+9x^3+4x^2$.
\end{lema}

\nt{\bf Proof.}
The proof of parts $(i)$ and $(ii)$ follow from   Theorems~\ref{theorem4} and ~\ref{theorem1} for vertex $u$ in graphs $Q_n^\bigtriangleup$ and $Q_n(2)$, respectively. Note that in these cases $p_u(G,x)=0$.

\nt $(iii)$ We use Theorems~\ref{theorem4} and ~\ref{theorem1} for vertex $u$ to obtain domination polynomial of $Q_n'$, then we have
\begin{eqnarray*}
D(Q_n',x)&=&(1+x)D(Q_n+e,x)+ x^2D(Q_{n-1}',x)-(1+x)xD(Q_{n-1}',x)\\
&=&(1+x)D(Q_n+e,x)- x^2D(Q_{n-1}',x).
\end{eqnarray*}

\nt $(iv)$ We use Theorems~\ref{theorem4} and ~\ref{theorem1} for vertex $w$ to obtain domination polynomial of $Q_n+e$, as shown  in figure~\ref{figure2} then we have
$D(Q_n+e,x)=xD((Q_n+e)/w,x)+ xD(Q_{n-1}',x)+xD(Q_{n-1},x).$ Now consider the graph $(Q_n+e)/w$ as shown in figure~\ref{figure2}. We use Theorems~\ref{theorem4} and ~\ref{theorem3} for $e=\{u, v\}$ to obtain $D((Q_n+e)/w,x)$, then we have
\begin{eqnarray*}
D((Q_n+e)/w,x)&=&D(Q_n,x) + \frac{x}{x-1}[D(Q_{n-1}^\bigtriangleup, x)+ D(Q_{n-1}^\bigtriangleup, x)- (Q_{n-1}^\bigtriangleup, x)\\
&&- D(Q_{n-1}^\bigtriangleup, x)- D(Q_{n-2}', x)- D(Q_{n-2}', x)+ xD(Q_{n-2}', x)+ xD(Q_{n-2}', x)]\\
&=& D(Q_{n}, x)+ 2xD(Q_{n-2}', x).
\end{eqnarray*}
\nt  Therefore the result follows.\quad\qed


\begin{teorem}\label{precurrence}
The domination polynomial of para-chain  $Q_n$  is given by
\begin{eqnarray*}
D(Q_n,x)&=& (x^3+2x^2+x)D(Q_{n-1},x)+ (x^3+2x^2)D(Q_{n-2},x)\\
&&+ (x^3+3x^2)D(Q_{n-2}',x)+ (2x^4+4x^3)D(Q_{n-3}',x),
\end{eqnarray*}
 with initial conditions $D(Q_1,x)=x^4+4x^3+6x^2$ and $D(Q_2,x)=x^7+7x^6+21x^5+29x^4+15x^3$.
\end{teorem}
\nt{\bf Proof.}
Consider the labeled $Q_n$ as shown in Figure \ref{pcactus}. We use Theorems~\ref{theorem4} and ~\ref{theorem1} for vertex $u_n$ to obtain the domination polynomial of $Q_n$. We have
\begin{eqnarray}\label{eq2}
D(Q_n,x)&=&x D(Q_{n-1}^\bigtriangleup, x)+  D(Q_{n-1}(2), x)+ x^2 D(Q_{n-2}', x)- (1+x)x D(Q_{n-2}', x)\nonumber\\
&=&x D(Q_{n-1}^\bigtriangleup, x)+  D(Q_{n-1}(2), x)- x D(Q_{n-2}', x).
\end{eqnarray}
\nt Therefore by parts $(i), (ii)$ and $(iv)$ of Lemma~\ref{lemma1} and equation (\ref{eq2}) we have
\begin{eqnarray*}
D(Q_n,x)&=&x ((1+x)D(Q_{n-1}+e,x) + xD(Q_{n-2}',x))+ x(D(Q_{n-1}+e,x)\\
&& + D(Q_{n-1},x) + D(Q_{n-2}',x))- x D(Q_{n-2}', x)\\
&=& (x^2+2x)D(Q_{n-1}+e,x) + x^2D(Q_{n-2}', x) + xD(Q_{n-1}, x)\\
&=& (x^2+2x)[x(D(Q_{n-1},x) + D(Q_{n-2},x)) + xD(Q_{n-2}',x)\\
&& + 2x^2D(Q_{n-3}',x)]+ x^2D(Q_{n-2}', x) + xD(Q_{n-1}, x)\\
&=& (x^3+2x^2+x)D(Q_{n-1},x)+ (x^3+2x^2)D(Q_{n-2},x)\\
&&+ (x^3+3x^2)D(Q_{n-2}',x)+ (2x^4+4x^3)D(Q_{n-3}',x).\quad\qed
\end{eqnarray*}

\subsection{Domination polynomial of ortho-chain square cactus graphs}

\begin{figure}[!h]
\hspace{3.9cm}
\includegraphics[width=8.3cm,height=2.7cm]{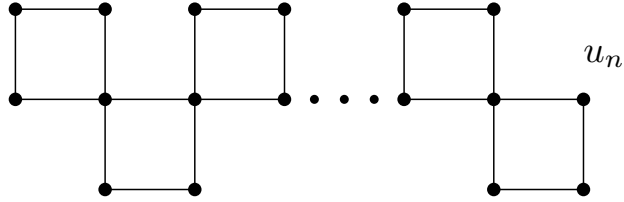}
\caption{ \label{scactus}  Labeled  ortho-chain square $O_n$. }
\end{figure}

\nt In this subsection we consider a ortho-chain of length $n$, $O_n$,  as shown in Figure \ref{scactus}. We shall
obtain a recurrence relation for the domination polynomial of $O_n$.

\nt We need the following Lemma for finding domination polynomial of the $O_n$.

\begin{figure}[!h]
\hglue2.5cm
\includegraphics[width=12cm,height=2cm]{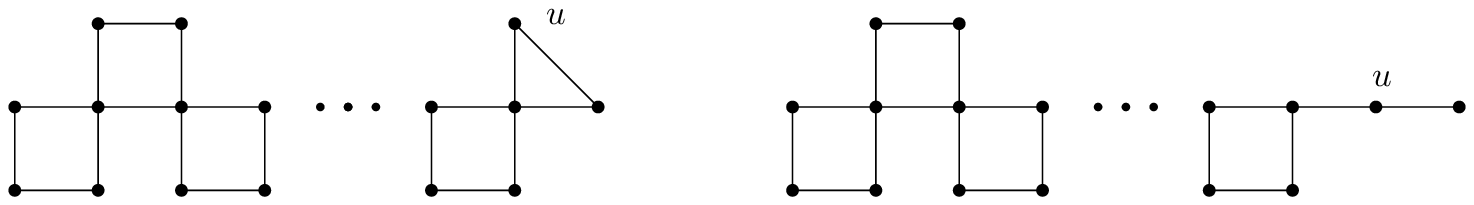}
\vglue5pt
\hglue2.5cm
\includegraphics[width=12cm,height=2cm]{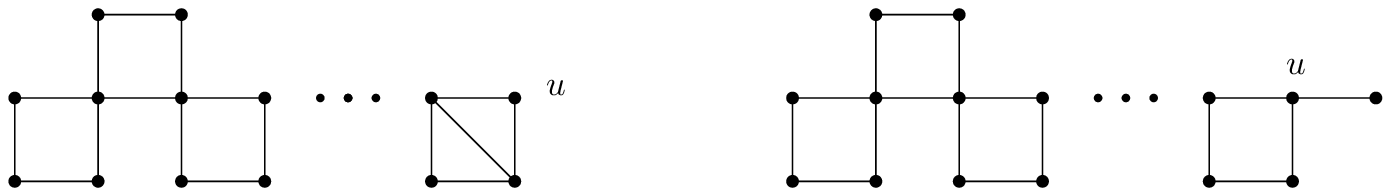}
\caption{\label{figure3}  Graphs $O_n^\bigtriangleup$, ~$O_n(2)$,  ~$O_n'$~and~ $O_n+e$, respectively.}
\end{figure}

\begin{lema}\label{lemma2}
\nt For graphs  in figure~\ref{figure3} we have:

\nt $(i)~ D(O_n^\bigtriangleup, x)=(1+x)D(O_n+e,x) + xD(O_{n-1}(2),x)$, where $D(O_0^\bigtriangleup, x)=x^3+3x^2+3x$.

\nt $(ii)~ D(O_n(2), x)=x(D(O_n+e,x) + D(O_n,x) + D(O_{n-1}(2),x))$, where $D(O_0(2), x)=x^3+3x^2+x$.

\nt $(iii)~ D(O_n', x)=(1+x)D(O_{n}^\bigtriangleup,x) - xD(O_{n-1}(2), x)$, where $D(O_0', x)= x^4+4x^3+6x^2+2x$.

\nt $(iv)~ D(O_n+e, x)=xD(O_n',x) +xD(O_{n-1}(2), x) + x^2D(O_{n-2}(2), x)$, where $D(O_1+e, x)=x^5+5x^4+9x^3+4x^2$.
\end{lema}

\nt{\bf Proof.}
The proof of parts $(i),~(ii)$ and $(iv)$ follow from   Theorems~\ref{theorem4} and ~\ref{theorem1} for vertex $u$ in graphs $O_n^\bigtriangleup,~O_n(2)$ and $O_n+e$,  respectively. Note that in these cases $p_u(G,x)=0$.

\nt $(iii)$ We use Theorems~\ref{theorem4} and ~\ref{theorem1} for $u$ in graphs $~O_n'$.
 Since $O_n'/u$ is isomorphic to $O_n' - u$ and $p_u(G,x)=xD(O_{n-1}(2),x)$. So we have the result.
\quad\qed


\begin{teorem}\label{srecurrence}
The domination polynomial of para-chain  $O_n$  is given by
$$D(O_n,x)=xD(O_{n-1},x)+ (x^2+2x)D(O_{n-1}+e,x)+x^2D(O_{n-2}(2),x),$$
with initial condition $D(O_1,x)=x^4+4x^3+6x^2.$
\end{teorem}
\nt{\bf Proof.}
Consider the labeled $O_n$ as shown in Figure \ref{scactus}. We use Theorems~\ref{theorem4} and ~\ref{theorem1} for vertex $u_n$
to obtain domination polynomial of $O_n$, then we have
\begin{eqnarray*}
D(O_n,x)&=&x D(O_{n-1}^\bigtriangleup, x)+  D(O_{n-1}(2), x)+ x^2 D(O_{n-2}(2), x)- (1+x)x D(O_{n-2}(2), x)\\
&=&x D(O_{n-1}^\bigtriangleup, x)+  D(O_{n-1}(2), x)- x D(O_{n-2}(2), x).
\end{eqnarray*}
\nt Therefore by parts $(i)$ and $(ii)$ of Lemma~\ref{lemma2} and this equation we have
\begin{eqnarray*}
D(O_n,x)&=&x ((1+x)D(O_{n-1}+e,x) + xD(O_{n-2}(2),x))+ x(D(O_{n-1}+e,x) \\
&&+ D(O_{n-1},x) + D(O_{n-2}(2),x))- x D(O_{n-2}(2), x)\\
&=& (x^2+2x)D(O_{n-1}+e,x) + x^2D(O_{n-2}(2), x) + xD(O_{n-1}, x).\quad\qed
\end{eqnarray*}


\begin{thebibliography}{99}



\bibitem{euro} S. Akbari, S. Alikhani and Y.H.  Peng, {\it Characterization of
graphs using domination polynomial}, Europ. J. Combin.,  Vol 31 (2010) 1714-1724.


\bibitem {saeid2}  S. Alikhani, {\it On the domination polynomials of non $P_4$-free graphs},
Iran. J. Math. Sci. Informatics, Vol. 8, No. 2 (2013) 49--55.

\bibitem{gcom} S. Alikhani, {\it The domination polynomial of a graph at $-1$}, Graphs  Combin., 29 (2013) 1175-1181.


\bibitem{saeid1}  S. Alikhani, Y.H. Peng, {\it Introduction to domination polynomial of a graph}, Ars Combin., to appear. Available at \texttt{http://arxiv.org/abs/0905.2251}.


\bibitem{doslicnew} K. Borissevich, T. Do\v{s}li\'{c}, Counting dominating sets in cactus chains, submitted. 

\bibitem{chellali} M. Chellali, {\it Bounds on the 2-domination number in cactus graphs}, Opuscula Math. 26
(2006) 5--12.

\bibitem{eslahchi} C. Eslahchi, E.S. Ansari, {\it A new protein domain assignment algorithm
based on the dominating set of a graph}, MATCH Commun. Math. Comput. Chem. 71 (2014) 445--456.

\bibitem{9} F. Harary, B. Uhlenbeck, {\it On the number of Husimi trees,} I, Proc. Nat. Acad. Sci. 39 (1953) 315--322.


\bibitem{domination} T.W. Haynes, S.T. Hedetniemi, P.J. Slater, {\it Fundamentals of domination in graphs}, Marcel Dekker, NewYork, 1998.





\bibitem{12} K. Husimi, {\it Note on Mayer's theory of cluster integrals}, J. Chem. Phys. 18 (1950) 682--684.



\bibitem {Kot} T. Kotek, J. Preen, F. Simon,P. Tittmann, M. Trinks, {\it Recurrence relations and splitting formulas for the domination polynomial},  Elec. J. Combin. 19(3) (2012), \# P47.


\bibitem{gcomkot} T. Kotek, J. Preen, P. Tittmann, {\it Subset-Sum representations of domination polynomials}, Graphs Combin.
DOI 10.1007/s00373-013-1286-z.


\bibitem{13} S. Majstorovi\'{c},  T. Do\v{s}li\'{c}, A. Klobu\v{c}ar, {\it $k$-domination on hexagonal cactus chains}, Kragujevac J. Math.,
Vol 36 No 2 (2012), Pages 335-–347.

\bibitem{mesi} R. E. Merrifield and H. E. Simmons, {\it Topological methods in Chemistry}, Wiley, New York, 1989.











\bibitem{prod} H. Prodinger and R. F. Tichy, {\it Fibonacci numbers of graphs}, Fibonacci Quart., 20 (1982) 16-21.

\bibitem{14} R. J. Riddell, {\it  Contributions to the theory of condensation}, Ph.D. Thesis, Univ. of Michigan, Ann Arbor, 1951.


\bibitem{dam} Z. Skupie\'{n}, {\it Majorization and the minimum number of dominating sets}, Discrete Appl. Math. 165 (2014) 295-302.


\bibitem{wag} S. Wagner, {\it Extemal trees with respect to Hosoya index and Merrifield-Simmons index}, MATCH Commun. Math.
Comput. Chem., 57 (2007) 221-233.


\bibitem{Wal} M. Walsh,  {\it The hub number of a graph}, Int. J. Math. Comput. Sci., 1 (2006) 117-124.




\end{thebibliography}
\end{document}